\documentclass{gtart_h}

\def\ifplaintex{\expandafter\ifx\csname documentclass\endcsname\relax}

\def\gtp{{\mathsurround=0pt\it $\cal G\mskip-2mu$eometry \&\ 
$\cal T\!\!$opology $\cal P\!$ublications}}  

\def\recd{{\small Received:\qua\receiveddate\ifx\reviseddate\relax
\else\qquad Revised:\qua\reviseddate\fi\par}} 

\def\lognumber#1{\def\thelognumber{#1}}
\def\volumenumber#1{\def\thevolumenumber{#1}}
\def\volumeyear#1{\def\thevolumeyear{#1}}
\def\papernumber#1{\def\thepapernumber{#1}}
\def\pagenumbers#1#2{\def\startpage{#1}\def\finishpage{#2}}
\def\published#1{\def\publishdate{#1}}

\def\received#1{\def\receiveddate{#1}}

\def\accepted#1{\def\accepteddate{#1}}

\def\asciiemail#1{\def\theasciiemail{#1}}
\def\asciiurl#1{\def\theasciiurl{#1}}

\long\def\asciiabstract#1{\long\def\theasciiabstract{#1}}

\let\\\par\let\thelognumber\relax\let\thevolumenumber\relax
\let\thepapernumber\relax\let\thevolumeyear\relax\let\startpage\relax
\let\finishpage\relax\let\publishdate\relax\let\receiveddate\relax
\let\reviseddate\relax\let\accepteddate\relax\let\theasciititle\relax
\let\theasciiauthors\relax
\let\theasciiabstract\relax

\let\theasciiemail\relax
\let\theasciiurl\relax

\ifplaintex
\font\logobig=cmssbx10 scaled 3836
\font\logomed=cmssbx10 scaled 2557
\else
\font\logobig=cmssbx10 scaled 4200
\font\logomed=cmssbx10 scaled 2800
\fi

\long\def\makeagttitle{   
\count0=\startpage
\agt\hfill      
\hbox to 45truept{\vbox to 0pt{\vglue -13truept{\logomed A\kern -.37em{\logobig 
T}\kern -.38em G}\vss}\hss}
\break
{\small Volume \thevolumenumber\ (\thevolumeyear)
\startpage--\finishpage\nl
Published: \publishdate}

\vglue .25truein

{\parskip=0pt\leftskip 0pt plus
1fil\def\\{\par\smallskip}{\Large\bf\thetitle}\par\medskip} \vglue
0.05truein

{\parskip=0pt\leftskip 0pt plus 1fil\def\\{\par}{\sc\theauthors}
\par\medskip}%
 
\vglue 0.03truein 

{\small\leftskip 25truept\rightskip 25truept{\bf Abstract}\stdspace\theabstract

{\bf AMS Classification}\stdspace\theprimaryclass
\ifx\thesecondaryclass\relax\else; \thesecondaryclass\fi\par
{\bf Keywords}\stdspace \thekeywords\par}\vglue 7truept

}   

\ifplaintex
\font\phead=cmsl9 scaled 950
\font\pnum=cmbx10 scaled 913
\font\pfoot=cmsl9 scaled 950
\headline{\vbox to 0pt{\vskip -4.5mm\line{\small\phead\ifnum
\count0=\startpage ISSN 1472-2739 (on-line) 1472-2747 (printed)
\hfill {\pnum\folio}\else\ifodd\count0\def\\{ }%
\ifx\theshorttitle\relax\thetitle\else\theshorttitle\fi\hfill{\pnum\folio}
\else\def\\{ and }{\pnum\folio}\hfill\ifx\theshortauthors\relax\theauthors
\else\theshortauthors\fi\fi\fi}\vss}}
\footline{\vbox to 0pt{\vglue 0mm\line{\small\pfoot\ifnum\count0=\startpage
\copyright\ \gtp\hfill\else
\agt, Volume \thevolumenumber\ (\thevolumeyear)\hfill\fi}\vss}}
\else
\font=cmsl9 scaled 1050
\font\lnum=cmbx10 
\font=cmsl9 scaled 1050
\makeatletter
\def\@oddhead{{\small\lhead\ifnum\count0=\startpage ISSN 1472-2739 
(on-line) 1472-2747 (printed)\hfill {\lnum\number\count0}\else\ifodd\count0
\def\\{ }\ifx\theshorttitle\relax \thetitle \else\theshorttitle\fi\hfill
{\lnum\number\count0}\else\def\\{ and }{\lnum\number\count0}
\hfill\ifx\theshortauthors\relax 
\theauthors\else\theshortauthors\fi\fi\fi}}\def\@evenhead{\@oddhead}
\def\@oddfoot{\small\lfoot\ifnum\count0=\startpage\copyright\ \gtp\hfill\else
\agt, Volume \thevolumenumber\ (\thevolumeyear)\hfill\fi}
\def\@evenfoot{\@oddfoot}
\makeatother
\fi
\let\maketitlepage\makeagttitle
\let\makeshorttitle\maketitlepage
\let\maketitle\maketitlepage

\newwrite\gtoutfile
\long\gdef\makeheadfile{  
{\def\\{, }\def\s{ }
\immediate\openout\gtoutfile head.xxx
\immediate\write\gtoutfile{Proxy-for: \ifx\theasciiauthors\relax
\theauthors\else\theasciiauthors\fi\s<\ifx\theasciiemail\relax\theemail\else\theasciiemail\fi>}
\immediate\write\gtoutfile{\noexpand\\}
\immediate\write\gtoutfile{Authors: \ifx\theasciiauthors\relax
\theauthors\else\theasciiauthors\fi}
{\def\\{ }\immediate\write\gtoutfile{Title: \ifx\theasciititle\relax
\thetitle\else\theasciititle\fi}}
\immediate\write\gtoutfile{Subj-class: GT or SG, GR etc}
\immediate\write\gtoutfile{MSC-class: \theprimaryclass\ifx\thesecondaryclass\relax\else, \thesecondaryclass\fi}
\immediate\write\gtoutfile{Journal-ref: Algebr. Geom. Topol. \thevolumenumber\s
(\thevolumeyear) \startpage-\finishpage}
\immediate\write\gtoutfile{Comments: Published by Algebraic and
Geometric Topology at}
\immediate\write\gtoutfile{\s\s\s  http://www.maths.warwick.ac.uk/agt/AGTVol\thevolumenumber/agt-\thevolumenumber-\thepapernumber.abs.html}
\immediate\write\gtoutfile{\noexpand\\}
\immediate\write\gtoutfile{}
\ifx\theasciiabstract\relax
\immediate\write\gtoutfile{\theabstract}\else
\immediate\write\gtoutfile{\theasciiabstract}\fi
\immediate\write\gtoutfile{}
\immediate\write\gtoutfile{\noexpand\\}
\immediate\write\gtoutfile{}
\immediate\closeout\gtoutfile}}  

\def\maketitlepage{\makeagttitle\makeheadfile}
\let\makeshorttitle\maketitlepage
\let\maketitle\maketitlepage

\lognumber{29}
\volumenumber{4}
\volumeyear{2004}
\papernumber{29}
\received{6 February 2004} 
\pagenumbers{623}{645}
\published{23 August 2004}
\accepted{21 August 2004}

\usepackage{amssymb,amsmath,amsxtra}
\usepackage[all]{xy}
\usepackage{mathrsfs}

\newtheorem{theorem}{Theorem}[section]
\newtheorem{cor}[theorem]{Corollary}
\newtheorem{lem}[theorem]{Lemma}
\newtheorem{prop}[theorem]{Proposition}

\newtheorem{con}[theorem]{Conjecture}
\newtheorem{prob}[theorem]{Problem}

\theoremstyle{definition}
\newtheorem{defi}[theorem]{Definition}
\newtheorem{example}[theorem]{Example}
\newtheorem{rem}[theorem]{Remark}

\numberwithin{equation}{section}

\def\ds{\displaystyle}
\def\:{\colon\thinspace}
\def\.{\cdot}
\def\o{\circ}
\def\<{\left\langle}
\def\>{\right\rangle}
\def\({\left(}
\def\){\right)}

\def\epsilon{\varepsilon}
\def\phi{\varphi}

\def\geq{\geqslant}

\def\lra{\longrightarrow}
\def\Lra{\Longrightarrow}
\def\ra{\rightarrow}
\def\bar#1{\overline{#1}}
\def\hat#1{\widehat{#1}}
\def\tilde#1{\widetilde{#1}}
\def\iso{\cong}

\def\ideal{\triangleleft}

\DeclareMathOperator{\F}{F}
\DeclareMathOperator{\id}{id}
\def\op{\mathrm{op}}

\DeclareMathOperator{\Ker}{Ker}\renewcommand{\ker}{\Ker}
\DeclareMathOperator{\End}{End}
\DeclareMathOperator{\Ext}{Ext}
\DeclareMathOperator{\Hom}{Hom}

\DeclareMathOperator{\RHom}{RHom}

\DeclareMathOperator{\HH}{HH}
\DeclareMathOperator{\THH}{THH}

\def\LMod#1{{}_{#1}\mathscr{M}}
\def\RMod#1{\mathscr{M}_{#1}}
\def\LDMod#1{{}_{#1}\mathscr{D}}
\def\RDMod#1{\mathscr{D}_{#1}}
\def\Z{\mathbb{Z}}
\def\Times_#1{{\ds\mathop{\times}_{#1}}}
\def\oTimes_#1{{\ds\mathop{\otimes}_{#1}}}
\def\Smash_#1{{\ds\mathop{\wedge}_{#1}}}

\def\k{\Bbbk}
\def\functorto{\lra}

\begin{document}
\title{Topological Hochschild cohomology and\\generalized Morita equivalence}
\authors{Andrew Baker\\Andrey Lazarev}
\addresses{Mathematics Department, Glasgow University, Glasgow G12 8QW,
Scotland\\Mathematics Department, Bristol University, Bristol, BS8 1TW,
England.}
\gtemail{\mailto{a.baker@maths.gla.ac.uk}\qua {\rm and}\qua
\mailto{a.lazarev@bristol.ac.uk}}
\asciiemail{a.baker@maths.gla.ac.uk, a.lazarev@bristol.ac.uk}
\gturl{\url{http://www.maths.gla.ac.uk/~ajb/}{\rm\qua 
and\qua}\url{http://www2.maths.bris.ac.uk/~maxal/}}
\asciiurl{http://www.maths.gla.ac.uk/ ajb/,
http://www2.maths.bris.ac.uk/ maxal/}
\begin{abstract}
We explore two constructions in homotopy category with algebraic precursors
in the theory of noncommutative rings and homological algebra, namely the
Hochschild cohomology of ring spectra and Morita theory. The present paper
provides an extension of the algebraic theory to include the case when $M$
is not necessarily a progenerator. Our approach is complementary to recent
work of Dwyer and Greenlees and of Schwede and Shipley.

A central notion of noncommutative ring theory related to Morita equivalence
is that of \emph{central separable} or \emph{Azumaya} algebras. For such an
Azumaya algebra $A$, its Hochschild cohomology $\HH^*(A,A)$ is concentrated
in degree~$0$ and is equal to the center of $A$. We introduce a notion of
\emph{topological Azumaya algebra} and show that in the case when the ground
$\mathbb{S}$-algebra $R$ is an Eilenberg-Mac~Lane spectrum of a commutative
ring this notion specializes to classical Azumaya algebras. A canonical example
of a topological Azumaya $R$-algebra is the endomorphism $R$-algebra $\F_R(M,M)$
of a finite cell $R$-module. We show that the spectrum of mod~$2$ topological
$K$-theory $KU/2$ is a nontrivial topological Azumaya algebra over the $2$-adic
completion of the $K$-theory spectrum $\widehat{KU}_2$. This leads to the
determination of $\THH(KU/2,KU/2)$, the topological Hochschild cohomology of
$KU/2$. As far as we know this is the first calculation of $\THH(A,A)$ for a
\emph{noncommutative} $\mathbb{S}$-algebra $A$.
\end{abstract}

\asciiabstract{%
We explore two constructions in homotopy category with algebraic precursors
in the theory of noncommutative rings and homological algebra, namely the
Hochschild cohomology of ring spectra and Morita theory. The present paper
provides an extension of the algebraic theory to include the case when $M$
is not necessarily a progenerator. Our approach is complementary to recent
work of Dwyer and Greenlees and of Schwede and Shipley.
A central notion of noncommutative ring theory related to Morita
equivalence is that of central separable or Azumaya algebras. For such
an Azumaya algebra A, its Hochschild cohomology HH^*(A,A) is
concentrated in degree 0 and is equal to the center of A.  We
introduce a notion of topological Azumaya algebra and show that in the
case when the ground S-algebra R is an Eilenberg-Mac Lane spectrum of
a commutative ring this notion specializes to classical Azumaya
algebras.  A canonical example of a topological Azumaya R-algebra is
the endomorphism R-algebra F_R(M,M) of a finite cell R-module.  We
show that the spectrum of mod 2 topological K-theory KU/2 is a
nontrivial topological Azumaya algebra over the 2-adic completion of
the K-theory spectrum widehat{KU}_2. This leads to the determination
of THH(KU/2,KU/2), the topological Hochschild cohomology of KU/2.  As
far as we know this is the first calculation of THH(A,A) for a
noncommutative S-algebra A.}

\maketitle

\section*{Introduction}
\addcontentsline{toc}{section}{Introduction}

The recent development of models for stable homotopy category with strictly
symmetric monoidal structure has led to the whole new outlook on the relationship
between algebra and homotopy theory. In this paper we explore two constructions
in homotopy category with algebraic precursors in the theory of noncommutative
rings and homological algebra, namely the Hochschild cohomology of ring spectra
and Morita theory.

Our approach to Morita theory was inspired by the paper of Dwyer
and Greenlees~\cite{DG}. Recall that the classical Morita theorem asserts that
the categories of $R$-modules and $\End_R(M)$-modules are equivalent, where~$R$
is a ring and $M$ is a faithfully projective $R$-module (or, in another
terminology, an $R$-progenerator, see~\cite{Bass}). Both the present paper
and~\cite{DG} provide an extension of this theorem to include the case when~$M$
is not necessarily a progenerator. Our approach is in some sense complementary
to that of Dwyer and Greenlees. Also closely related to our work is that of
Schwede and Shipley~\cite{Schwede-Shipley}, in which Morita theory for modules
over `ring spectra with several objects' is developed.

One of the central notions of noncommutative ring theory closely related to
Morita equivalence is that of a \emph{central separable} algebra or
\emph{Azumaya} algebra,~\cite{AG}. A characteristic property of an Azumaya
algebra~$A$ is that its Hochschild cohomology $\HH^*(A,A)$ is concentrated
in degree~$0$ and is equal to the center of $A$. We introduce the notion of
a \emph{topological Azumaya algebra} and show that in the case when the ground
$\mathbb{S}$-algebra $R$ is an Eilenberg-Mac~Lane spectrum of a commutative
ring such algebras specialize to classical Azumaya algebras.

A canonical example of a topological Azumaya $R$-algebra is the endomorphism
$R$-algebra $\F_R(M,M)$ of a finite cell $R$-module. By analogy with the
theory of the Brauer group, it is natural to consider such Azumaya $R$-algebras
\emph{trivial}. We show that $KU/2$, the spectrum of topological $K$-theory
mod~$2$ is a nontrivial topological Azumaya algebra over $\widehat{KU}_2$, the
topological $K$-theory spectrum completed at~$2$.

This leads to the determination of $\THH(KU/2,KU/2)$, the topological Hochschild
cohomology of $KU/2$. As far as we know this is the first calculation of $\THH(A,A)$
for a \emph{noncommutative} $\mathbb{S}$-algebra $A$. A somewhat surprising outcome
is that
\[
\THH(KU/2,KU/2)\simeq\widehat{KU}_2.
\]

\bigskip
We would like to thank John Greenlees and Stefan Schwede for helpful comments.

\subsection*{Notation and conventions}
\addcontentsline{toc}{subsection}{Notation and conventions}

We choose to work throughout in the category of $\mathbb{S}$-modules of~\cite{EKMM},
where $\mathbb{S}$ denotes the sphere spectrum. But we note that our constructions
should be applicable in any version of the underlying model category for spectra
with symmetric monoidal structure such as the category of symmetric spectra~\cite{SymmSpect}.
We will use the symbols $\iso$ and $\simeq$ to denote point-set isomorphism and
weak equivalence, respectively. We also use the term \emph{$q$-cofibrant} in the
sense of~\cite{EKMM}.

If $R$ is an $\mathbb{S}$-algebra, we will often work in the category of left/right
$R$-modules $\LMod{R}$/$\RMod{R}$ and their derived homotopy categories
$\LDMod{R}$/$\RDMod{R}$ as described in~\cite{EKMM}. To ease notation, we will
usually write $M\wedge N$ for the $R$-smash product $M\wedge_R N$, $\F_R(M,N)$
or $\F(M,N)$ for the function object $\F_{\LMod{R}}(M,N)$ or $\F_{\RMod{R}}(M,N)$.
We will also write $[M,N]_*$ for $[M,N]^R_*=\pi_*\F_R(M,N)$ when it is clear
that no ambiguity will result.

We will denote the Spanier-Whitehead dual of an $R$-module $M$ by $M^*$,
so $M^*=\F_R(M,R)$.

Following~\cite{D&P,LMS,HPS}, we call an $R$-module $M$ \emph{dualizable}
if it is weakly equivalent to a retract of a finite cell $R$-module. Our
justification for this terminology lies in the fact that for objects of
$\mathscr{D}_R$, the condition of being strongly dualizable in the sense
of~\cite{D&P,LMS} is equivalent to being a retract of a finite cell $R$-module;
see~\cite[\S2]{HPS},~\cite[proposition~2.1]{FLM} and~\cite{JPM:EqtHtpyCohThy}.
It has been pointed out by the referee that a better word to use might have
been `finite'.

Suppose that $A$ is either an $R$-algebra or an $R$-ring spectrum for a
commutative $\mathbb{S}$-algebra $R$. Write $\phi\:A\wedge A\lra A$ for
the product and $\tau\:A\wedge A\lra A\wedge A$ for the switch map. Then we
also have the opposite $R$-algebra $A^\op$ which consists of the $R$-module~$A$
with the product map $\phi^\op=\phi\circ\tau\:A\wedge A\lra A$. By considering~$A$
as a left $A\wedge A^\op$-module or a left $A\wedge A^\op$-module spectrum,
we have the action map in $\LMod{R}$
\[
\phi\circ(\id\wedge\phi\circ\tau)\:A\wedge A^\op\wedge A\lra A
\]
which has an adjoint
\[
\mu\:A\wedge A^\op\lra\F(A,A).
\]
This map $\mu$ passes to an adjoint in $\LDMod{R}$ which we also denote by~$\mu$.
It is straightforward to check that the following is true.
\begin{lem}\label{lem:ActionMap}
The action map $\mu\:A\wedge A^\op\lra\F(A,A)$ is
\begin{itemize}
\item
a map of $R$-algebras if $A$ is an $R$-algebra;
\item
a map of  $R$-ring spectra if $A$ is an $R$-ring spectrum.
\end{itemize}
\end{lem}

Let $\mathscr{C}$ and $\mathscr{D}$ be categories. Given two functors
$F\:\mathscr{C}\functorto\mathscr{D}$ and
$G\:\mathscr{D}\functorto\mathscr{C}$ we say that $F$ and $G$ form
an \emph{adjoint pair} if $F$ is left adjoint to $G$, or equivalently if
$G$ is right adjoint to $F$. We write $\mathscr{C}\simeq\mathscr{D}$ if
there is an equivalence $\mathscr{C}\functorto\mathscr{D}$.

\section{Generalized Morita theory}\label{sec:GenMorThy}

In this section we describe a rather general form of Morita theory. Our
approach should be compared with that of Dwyer and Greenlees~\cite{DG}. Here
we assume that $R$ is a not necessarily commutative $\mathbb{S}$-algebra
and $E$ is a left $R$-module. We have the $\mathbb{S}$-algebra $A=\F_R(E,E)$
for which $E$ is a left $A$-module. We will make use of the notions in the
following definition, the first of which is perhaps more familiar.
\begin{defi}\label{defi:Colocal}
Let $E$ be a left $R$-module.
\begin{itemize}
\item
A right $R$-module $M$ is \emph{local with respect to $E$} or \emph{$E$-local},
if for any $E$-acyclic right $R$-module $X$, $\LDMod{R}(X,M)$ is trivial.
\item
A right $R$-module $N$ is \emph{colocal with respect to $E$} or \emph{$E$-colocal},
if for any $E$-acyclic right $R$-module $Y$, $\LDMod{R}(N,Y)$ is trivial.
\end{itemize}
\end{defi}

In the following we will work with both of the categories of left $A$-modules
$\LMod{A}$ and right $R$-modules $\RMod{R}$. Consider the following functors
\begin{align*}
F_0\:\RMod{R}\functorto\LMod{A};&\quad X\longmapsto
X\wedge_RE, \\ G_0\:\LMod{A}\functorto\RMod{R};&\quad
Y\longmapsto\F_A(E,Y) \\
\end{align*} It is straightforward to see that $F_0$ and $G_0$
form an adjoint pair. However, although $F_0$ induces a functor
$\RDMod{R}\functorto\LDMod{A}$, it is not clear that $G_0$ induces
a functor $\LDMod{A}\functorto\RDMod{R}$. To get around this
issue, we choose a cofibrant replacement $\tilde E$ for $E$ in
$\LMod{A}$ and introduce the functors \begin{align*}
F\:\RMod{R}\functorto\LMod{A};&\quad X\longmapsto X\wedge_R\tilde
E, \\ G\:\LMod{A}\functorto\RMod{R}; &\quad
Y\longmapsto\F_A(\tilde E,Y). \end{align*}
\begin{theorem}\label{thm:mor}  The functors $F$ and $G$ form an
adjoint pair. Furthermore, they induce an adjoint pair of functors
$\tilde F\:\RDMod{R}\functorto\LDMod{A}$ and $\tilde
G\:\LDMod{A}\functorto\RDMod{R}$.

 If $E$ is a dualizable and $q$-cofibrant left $R$-module,
then $\tilde F$ and $\tilde G$ form an adjoint pair of
equivalences between the homotopy categories of right $E$-local
$R$-modules and left $A$-modules.

\end{theorem}
\begin{proof}
The first two statements are immediate. For the remaining part we first need
to show that the image of the functor $\tilde G$ lies within the subcategory
of $E$-local objects in $\RDMod{R}$. Let $M$ be an $E$-acyclic right $R$-module,
which means that $M\wedge_R E$ is null in $\LDMod{S}$ and hence in $\LDMod{A}$
since $\LDMod{A}(A,M\wedge_R E)\iso\LDMod{S}(S,M\wedge_R E)$. We may as well
assume that $M$ is $q$-cofibrant in $\RMod{R}$. Then
\[
\RDMod{R}(M,\tilde G(Y))=\RDMod{R}(M,\F_A(\tilde
E,Y))=\LDMod{A}(M\wedge_R\tilde E,Y) =\LDMod{A}(M\wedge_RE,Y)=0.
\]

We will show that $\tilde F\circ\tilde G\simeq\id_{\LDMod{A}}$ when
$E$ is dualizable. Then  $E^*=\F_{\LMod{R}}(E,R)$ is a right $R$-module
for which we choose a cofibrant replacement $W$. Now
$\F_{\RMod{R}}(W,R)\simeq\tilde E$ in $\LMod{R}$, since $E$ is dualizable
therein. Then in $\LMod{A}$ we have
\begin{align*}
F\circ G(Y)=\F_{\LMod{A}}(\tilde E,Y)\wedge_R\tilde E
&\simeq\F_{\LMod{A}}(\tilde E,Y)\wedge_R\F_{\RMod{R}}(W,R) \\
&\simeq\F_{\RMod{R}}(W,\F_{\LMod{A}}(\tilde E,Y)) \\
&\simeq\F_{\LMod{A}}(W\wedge_R\tilde E,Y) \\
&\simeq\F_{\LMod{A}}(W\wedge_RE,Y) \\
&\simeq\F_{\LMod{A}}(\F_{\LMod{R}}(E,E),Y) \\ &\simeq Y.
\end{align*}
Here we have used the condition that $E$ is $q$-cofibrant in $\LMod{R}$
and the evident equivalences of left $A$-modules
\[
A=\F_{\LMod{R}}(E,E)\simeq\F_{\LMod{R}}(E,R)\wedge_{R}E \simeq
W\wedge_{R}\tilde E.
\]

Conversely, let $X$ be an $E$-local right $R$-module. To prove
that $\tilde G\circ\tilde F(X)\simeq X$ it is enough to show that
$\tilde G\circ\tilde F(X)\wedge_R E\simeq X\wedge_R E$. But this
follows from the previously established fact that $\tilde
F\circ\tilde G\cong\id_{\LDMod{A}}$. \end{proof}
\begin{rem}\label{rem:localization}
Note that in the case when $E$ is a dualizable $q$-cofibrant left
$R$-module, the functor $\tilde{G}\circ\tilde{F}$ considered as an
endofunctor on $\mathscr{D}_R$ is isomorphic to the $E$-localization
functor $X\mapsto X_E$ associating to $X\in\mathscr{D}_R$ a left
$R$-module $X_E$ which is $E$-local and $E$-equivalent to $X$. Indeed,
$\tilde{G}\circ\tilde{F}$ is $E$-local since the image of $\tilde{G}$
is $E$-local and the arguments at the end of the proof of Theorem~\ref{thm:mor}
show that the natural map $X\lra\tilde{G}\o\tilde{F}(X)$ is an $E$-equivalence.
\end{rem}

Now consider the following pair of functors:
\begin{align*}
H\:\LMod{A}\functorto\RMod{R};&\quad
Y\longmapsto\F_{\LMod{R}}(E,R)\wedge_AY, \\
F\:\RMod{R}\functorto\LMod{A};&\quad X\longmapsto X\wedge_R E.
\end{align*}
The functors $H$ and $G$ are \emph{not} adjoint on the point-set level,
however they are becoming such after passing to the derived categories.
More precisely we have the following proposition whose proof is analogous
to that of Theorem~\ref{thm:mor} and is left to the reader.
\begin{prop}\label{prop:add}
Let $E$ be a dualizable left $R$-module. Then the functors $H$ and $F$
induce functors $\tilde H\:\LDMod{A}\functorto\RDMod{R}$ and
$\tilde F\:\RDMod{R}\functorto\LDMod{A}$ which form an adjoint pair of
equivalences between the homotopy categories of left $A$-modules and
right $E$-colocal $R$-modules.
\end{prop}
\begin{rem}\label{rem:add}
Combining Theorem~\ref{thm:mor} with Proposition~\ref{prop:add}, we find
that the homotopy categories of right $E$-local and right $E$-colocal
$R$-modules are equivalent provided $E$ that is a dualizable right
$R$-module. This observation also appears in~\cite[theorem~3.3.5(g)]{HPS}.
\end{rem}

Our Theorem~\ref{thm:mor} is strikingly similar in formulation and proof
to theorem~2.1 of Dwyer and Greenlees~\cite{DG}, however it is essentially
different. In particular, in their context the appropriate analogues of
$E$-local and $E$-acyclic modules are those of \emph{$E$-torsion} and
\emph{$E$-trivial} modules given in the following definition.
\begin{defi}\label{defi:E-trivial/E-torsion}
A left $R$-module $N$ is \emph{$E$-trivial} if $\LDMod{R}(E,N)$ is trivial.
A left $R$-module $M$ is \emph{$E$-torsion} if $\LDMod{R}(M,N)$ is trivial
for all $E$-trivial $N$.
\end{defi}

For completeness, we reproduce part of~\cite[theorem~2.1]{DG}, even though
in what follows we stick to the version of Theorem~\ref{thm:mor}.

The notions of $E$-trivial and $E$-torsion left $R$-modules are in some ways
dual to $E$-acyclic and $E$-local right $R$-modules. We will make use of the
functors
\begin{align*}
T\:\RMod{A}\functorto\LMod{R};&\quad X\longmapsto X\wedge_A E, \\
S\:\LMod{R}\functorto\RMod{A};&\quad M\longmapsto\F(E,M).
\end{align*}

Now we can formulate a theorem which is `dual' to Theorem~\ref{thm:mor}
and appears in~\cite{DG}.
\begin{theorem}\label{thm:S-T-Adjoint}
The functors $T$ and $S$ form an adjoint pair. Furthermore, they induce
an adjoint pair of functors $\tilde T\:\LDMod{A}\functorto\LDMod{R}$ and
$\tilde S\:\LDMod{R}\functorto\RDMod{A}$.

If $E$ is a dualizable and $q$-cofibrant left $R$-module, then $\tilde T$ and
$\tilde S$ form an adjoint pair of equivalences between the homotopy categories
of left $E$-torsion $R$-modules and right $A$-modules.
\end{theorem}
\begin{rem}\label{rem:S-T-Adjoint}
The remaining part of the Dwyer-Greenlees theorem is `dual' to Proposition~\ref{prop:add}.
\end{rem}

We now return to Theorem~\ref{thm:mor}. The main advantage of this
over the classical Morita theorem (apart from the fact that we are
working in the topological framework) is that $E$ is not required
to be a generator in $\LDMod{R}$. We cannot hope that the
categories $\RDMod{R}$ and $\LDMod{A}$ are equivalent if $E$ is
not dualizable. However for certain right $R$-modules $M$ it is
still possible that the unit of the adjunction
\[
M\longmapsto\tilde G\tilde F(M)=\F_A(E,M\wedge_RE)
\]
is an $E$-equivalence (without the assumption that $E$ is
dualizable). We are especially interested when this is the case
for $M=R$.
\begin{defi}\label{defi:E-completion}
For a right $R$-module $M$ define its $E$-completion $E\sphat M$
to be the right $R$-module $\tilde G\tilde F(M)$ in $\LDMod{R}$.
\end{defi}

We now make a few easy but important observations which sometimes
help to decide whether the $E$-completion agrees with
$E$-localization. Since the functor $\tilde G$ takes its values in
the subcategory of $E$-local right $R$-modules the $E$-completion
of any right $R$-module will be $E$-local. The unit and counit
adjunction morphisms $?\longmapsto\tilde G\tilde F(?)$ and $\tilde
F\tilde G(?)\longmapsto?$ determine the pair of maps
\[
E\wedge_RM=FM\lra\tilde F\tilde G\tilde FM=E\wedge_R E\sphat M\lra
FM=E\wedge_RM,
\]
whose composite is obviously the identity map. Therefore the
$\mathbb{S}$-module $E\wedge_RM$ is a retract of $E\wedge_R E\sphat M$.
In particular, the $E$-homology of $M$ is a direct summand of the
$E$-homology of $E\sphat M$. Therefore a map $M\lra N$ is an
$E$-equivalence if and only if the induced map $E\sphat M\lra E\sphat N$
is an equivalence.

Now observe that for a right $R$-module $M$, the canonical map
$M\lra E\sphat M$ is an $E$-equivalence if and only if the map
$E_*M\lra E_*E\sphat M$ is an isomorphism, which happens exactly when
$E\sphat E\sphat M=E\sphat M$. In other words, the $E$-completion agrees
with $E$-localization in those cases when the $E$-completion is idempotent.
The following two purely algebraic examples are instructive.
\begin{example}\label{ex:1}
Let $R=\k[y]$ be the polynomial algebra over a field $\k$ on one generator
$y$ in degree~$0$ and $E=\k$ considered as a left $R$-module. Then
$A=\RHom_R(\k,\k)$ is easily seen to be $\Lambda(x)$, the exterior algebra
on one generator in degree~$-1$. Indeed, the unit map $\k\lra\k[y]$ gives
rise to an augmentation $A\lra k$ and it is easy to see that any
$\k$-augmented differential graded $\k$-algebra whose homology algebra is
$\Lambda(x)$ is formal, i.e.\ is quasi-isomorphic to its homology. Since~$E$
is a dualizable $R$-module we see that the category of $E$-local $R$-modules
is equivalent to the category of $A$-modules. The canonical map $M\lra E\sphat M$
is the $E$-localization by Remark~\ref{rem:localization}. Letting $M=R$,
we see that
\[
R_E= E\sphat R=\RHom_{\Lambda(x)}(\k,\k)=\k[[y]].
\]
\end{example}
\begin{example}\label{ex:2}
Let $R=\Lambda_\k(x)$ to be the exterior algebra over a field~$\k$ on
one generator of degree~$-1$ and $E=\k$ considered as a left $R$-module.
Then $A=\RHom_R(\k,\k)$ is easily seen to be the power series algebra
over $\k$ on one generator of degree~$0$. Since $\k$ is \emph{not} a
dualizable $R$-module (it has infinite homological dimension) we cannot
conclude that the category of left $R$-modules is equivalent to the
category of right $A$-modules. However we can analyze the situation
using Theorem~\ref{thm:S-T-Adjoint}. It is easy to see that the algebra
$R=\Lambda_\k(x)$ is quasi-isomorphic to $\RHom_A(\k,\k)$. Clearly the
(left) $A$-module $\k$ is dualizable. Therefore the category of $R$-modules
is equivalent to the category of torsion $A$-modules. In particular the
canonical map $R\lra E\sphat R$ is an equivalence even though $E$ is not
a dualizable $R$-module.
\end{example}

Here are some topological examples suggestive of double centralizer theorems
of algebra which we briefly recall, for details see~\cite{Lang,Pierce}.

Let $A$ be a $\k$-algebra over some field $\k$ and let $M$ be a simple
right $A$-module. Then $B=\End_AM$ acts naturally on the left of $M$
there is a canonical monomorphism $A\lra B$ which allows us to view $A$
as a $\k$-subalgebra of $B$. Then $\End_BM=A$.

In the topological context there are similar phenomena, at least in the
homotopy category. Such examples were hinted at by the first author
in~\cite{AB:Ainfty}, although the details required further elaboration of
the theory of $\mathbb{S}$-algebras and their modules to be made rigorous;
they are also related to calculations in~\cite{BaL}. For related ideas also
see~\cite{DGI}.
\begin{example}\label{ex:En}
From~\cite{GH}, it is known that for each prime~$p$ there is a commutative
$\mathbb{S}$-algebra $E_n$ whose homotopy ring is
\[
(E_n)_*=\pi_*E_n=
\mathbb{W}(\mathbb{F}_{p^n})[[u_1,\ldots,u_{n-1}]][u_n,u_n^{-1}],
\]
where $|u_1|=\cdots=|u_{n-1}|=0$, $|u_n|=2$ and $\mathbb{W}(\mathbb{F}_{p^n})$
is the ring of Witt vectors of the finite field $\mathbb{F}_{p^n}$. There is
also a $E_n$-ring spectrum $K_n$ with
\[
(K_n)_*=(E_n)_*/(u_1,\ldots,u_{n-1})=\mathbb{F}_{p^n}[u_n,u_n^{-1}].
\]
Of course, $K_n$ is s finite cell $E_n$-module. A straightforward calculation
with the hypercohomology spectral sequence of~\cite{EKMM} yields
\[
\pi_*\F_{E_n}(K_n,K_n)=\Lambda_{(K_n)_*}(Q^0,\ldots,Q^{n-1}),
\]
the exterior algebra over $(K_n)_*$ on Bockstein operations $Q^k$ of degree~$1$
associated with the elements $u_k$. The spectrum $A_n=\F_{E_n}(K_n,K_n)$ is itself
an $E_n$-algebra and $K_n$ is a left $A_n$-module. Another application of the
hypercohomology spectral sequence gives a spectral sequence
\[
\mathrm{E}_2^{*\,*}=\Ext_{(A_n)_*}^{s\,t}(K_n,K_n)\Lra\pi_*\F_{A_n}(K_n,K_n).
\]
This time we find that
\[
\mathrm{E}_2^{*\,*}=(K_n)_*[p_0,\ldots,p_{n-1}],
\]
where $p_k\in\mathrm{E}_2^{1\,0}$ represents the extension associated with $Q^k$.
Using Theorem~\ref{thm:mor} we see that $E_n\lra\F_{A_n}(K_n,K_n)$ is a weak
equivalence of $E_n$-ring spectra.
\end{example}
\begin{example}\label{ex:MU}
For a prime~$p$, the Eilenberg-Mac~Lane spectrum $H\mathbb{F}_p$ is naturally
a module over the commutative $\mathbb{S}$-algebra $MU\sphat_p$\,. Arguing in
similar fashion to Example~\ref{ex:En}, we find that for
$B_p=\F_{MU\sphat_p}(H\mathbb{F}_p,H\mathbb{F}_p)$,
\[
\pi_*B_p=\Lambda_{\mathbb{F}_p}(Q^k:k\geq0),
\]
where $Q^k$ is a Bockstein operation associated with a polynomial generator $x_k$
for $MU_*$ of degree $2k$ or $x_0=p$ when $k=0$. The natural map
$MU\sphat_p\lra\F_{B_p}(H\mathbb{F}_p,H\mathbb{F}_p)$ induces a weak equivalence
of $MU\sphat_p$-algebras. In this example we make use of the fact that for the
maximal ideal $I=(x_k:k\geq0)\ideal MU_*$\,,
\[
(MU_*)\sphat_I={MU_*}\sphat_p\,.
\]
\end{example}

More generally, using results of~\cite{BaL} it follows that for any ideal
$J\ideal MU_*$ generated by a regular sequence there is an $MU$-module $MU/J$ which
is a module over the $MU$-algebra $MU\sphat_J$. Then for $B_J=\F_{MU\sphat_J}(MU/J,MU/J)$,
the natural map
\[
MU\sphat_J\lra\F_{B_J}(MU/J,MU/J)
\]
is an equivalence of $MU\sphat_J$-ring spectra. There is also an analogue of this
example for the $\mathbb{S}$-algebra $BP$ and its maximal ideal $I_\infty=(u_k:k\geq0)$
for a collection of polynomial generators $u_k\in BP_{2(p^k-1)}$ and $v_0=p$. Under
the forgetful map
$\F_{BP\sphat_p}(H\mathbb{F}_p,H\mathbb{F}_p)\lra\F_S(H\mathbb{F}_p,H\mathbb{F}_p)$,
$\pi_*\F_{BP\sphat_p}(H\mathbb{F}_p,H\mathbb{F}_p)$ realizes the Bockstein subalgebra
of the usual mod~$p$ Steenrod algebra which can be identified with
$\pi_*\F_S(H\mathbb{F}_p,H\mathbb{F}_p)$.

\section{Azumaya $R$-algebras}\label{sec:AzumayaAlg}

In this section we introduce the notion of topological Azumaya $R$-algebras and
prove some of their properties. Azumaya algebras (also known as central separable
algebras) have been studied extensively in ring theory, especially in connection
with the Brauer group, see for example~\cite{AG}.

We assume that $R$ is a $q$-cofibrant commutative $\mathbb{S}$-algebra. Recall
that we usually write $\wedge$ in place of $\wedge_R$ and $\F$ in place $\F_R$.

Let $A$ be a $q$-cofibrant $R$-algebra and let $\tilde A$ be a $q$-cofibrant
replacement for it as an $A\wedge A^\op$-module. Then we can form the topological
Hochschild cohomology $R$-module of $A$ with coefficients in itself,
\[
\THH_R(A,A)=\F_{A\wedge A^\op}(\tilde A,\tilde A).
\]
Clearly $\THH_R(A,A)$ is an $R$-algebra with respect to the composition product.
\begin{defi}\label{defi:AzumayaAlg}
An $R$-algebra $A$ is a \emph{topological Azumaya algebra} provided the following
three conditions hold.
\begin{enumerate}
\item
$A$ is dualizable as a left $A\wedge A^\op$-module.
\item
$A\wedge A^\op$ is $A$-local as a left module over itself.
\item
$R\lra\THH_R(A,A)$ is a weak equivalence.
\end{enumerate}
\end{defi}
\begin{rem}\label{rem:AzumayaAlg}
Recall that an algebra $A$ over a commutative ring $R$ is an \emph{Azumaya
algebra} if and only if it is $R$-central (this corresponds to our condition~(3))
and a finitely generated (this corresponds to our condition~(1)) projective
generator in the category of left $A\otimes_R A^\op$-modules (this corresponds
to our condition~(2)).
\end{rem}

Our next result provides a convenient characterization of topological Azumaya
algebras. Recall the action map $\mu\:A\wedge A^\op\lra\F(A,A)$ of
Lemma~\ref{lem:ActionMap} which in this situation is a map of $R$-algebras.
\begin{prop}\label{prop:AzumayaAlg-Char}
An $R$-algebra $A$ is a topological Azumaya algebra if and only if the following
three conditions hold.
\begin{enumerate}
\item
$A$ is a dualizable $R$-module.
\item
$R$ is $A$-local as a module over itself.
\item
$\mu\:A\wedge A^\op\lra\F(A,A)$ is a weak equivalence.
\end{enumerate}
\end{prop}

\proof
Suppose that the three conditions (1),(2),(3) of Proposition~\ref{prop:AzumayaAlg-Char}
are satisfied. According to Theorem~\ref{thm:mor}, the categories
$\LDMod{A\wedge A^\op}\simeq\LDMod{\F(A,A)}$ and $\RDMod{R}$ are
equivalent via the functors $\tilde F$ and $\tilde G$ defined by
\[
\tilde{F}\:X\longmapsto X\wedge\tilde A,
\qquad
\tilde{G}\:Y\longmapsto\F_{A\wedge A^\op}(\tilde A,Y),
\]
for an $R$-module $X$ and a left $A\wedge A^\op$-module $Y$. The image
of $R$ under $\tilde F$ is the left $A\wedge A^\op$-module $\tilde A$
and it follows that $\tilde A$, and hence $A$, are both dualizable in
$\LMod{A\wedge A^\op}$ since $R$ is dualizable as module over itself.

Denote the $A$-localization of $R$ as an $R$-module by $R_A$ and that
of $A\wedge A^\op$ as an $A\wedge A^\op$-module by $(A\wedge A^\op)_A$.
Now using (3) and Theorem~\ref{thm:mor} we have
\[
\THH_R(A,A)\simeq\F_{\F(A,A)}(\tilde A,\tilde A)\simeq R_A\simeq R
\]
and
\[
A\wedge A^\op\simeq\F(A,A)\simeq\F_{\THH_R(A,A)}(\tilde A,\tilde A)
                                             \simeq(A\wedge A^\op)_A.
\]
This proves that $A$ is an topological Azumaya $R$-algebra.

Conversely, suppose that the conditions (1),(2),(3) of
Definition~\ref{defi:AzumayaAlg} are satisfied. By
Theorem~\ref{thm:mor}, the categories $\RDMod{A\wedge A^\op}$ and
$\LDMod{\THH_R(A,A)}\simeq\LDMod{R}$ are equivalent via the
functors $\tilde{F}$ and $\tilde{G}$ defined by
\[
\tilde{F}\:X\longmapsto X\wedge_{A\wedge A^\op}\tilde A,
\qquad
\tilde{G}\:Y\longmapsto\F(\tilde A,Y).
\]
Here $X$ is a right $A\wedge A^\op$-module and $Y$ is a left $R$-module.
The image of $A\wedge A^\op$ is the $R$-module $A$ and it follows that
$A$ is a dualizable $R$-module. Furthermore
\[
A\wedge A^\op\simeq(A\wedge A^\op)_A\simeq
                 \F_{\THH_R(A,A)}(\tilde A,\tilde A)\simeq\F(A,A).
\]
Finally we have
$$
R_A\simeq\F_{\F(A,A)}(\tilde A,\tilde A)\simeq
         \F_{A\wedge A^\op}(\tilde A,\tilde A)=\THH_R(A,A)\simeq R.
\eqno{\qed}
$$

\begin{rem}\label{rem:AzumayaAlg-Char}
Proposition~\ref{prop:AzumayaAlg-Char} is an analogue of the standard
characterization of an Azumaya algebra $A$ over a ring~$R$ as an
$R$-progenerator for which the natural map $A\otimes_R A^\op\lra\Hom_R(A,A)$
is an isomorphism.
\end{rem}
\begin{rem}\label{rem:AzumayaAlg-2}
Condition (3) of Proposition~\ref{prop:AzumayaAlg-Char}
should not be confused with the standard characterization of a dualizable
$R$-module $A$ by the weak equivalence $A\wedge A^*\simeq\F(A,A)$. We do
not know whether this condition actually implies dualizability. However,
dualizability itself certainly does not imply (3) as the case of a homotopy
commutative finite cell $R$-algebra makes clear.
\end{rem}

Now let $R=H\k$, the Eilenberg-Mac~Lane spectrum corresponding to
a commutative ring $\k$. Then the categories of $H\k$-modules and
$H\k$-algebras are Quillen equivalent to the categories of complexes
of $\k$-modules and differential graded $\k$-algebras respectively.
Recall that a \emph{perfect complex} of $\k$-modules is a finite
length complex consisting of finitely generated projective modules.
This is an analogue of a dualizable module in the algebraic category
of complexes of $\k$-modules. Under the above mentioned equivalence,
our definition of a generalized Azumaya algebra specializes to the
following.
\begin{defi}\label{defi:genAzumayaalg}
A \emph{generalized Azumaya $\k$-algebra} is a differential graded
$\k$-algebra $A$ such that
\begin{enumerate}
\item
$A$ is quasi-isomorphic to a perfect complex of $\k$-modules;
\item
$\k$ is $A$-local as a module over itself;
\item
the canonical map
$A\ds{\mathop{\otimes}^{\mathrm{L}}}_\k A^\op\lra\RHom_\k(A,A)$ is an
isomorphism in the derived category of $\k$-modules.
\end{enumerate}
\end{defi}
Recall that the classical Azumaya $\k$-algebra is defined as a
$\k$-algebra which is finitely generated and faithfully projective
$\k$-module $A$ such that the canonical map
$A\ds{\mathop{\otimes}^{\mathrm{L}}}_\k A^\op\lra\Hom_\k(A,A)$ is
an isomorphism. Then we have the following
\begin{prop}\label{prop:Azumaya}
Any Azumaya $\k$-algebra is a generalized Azumaya $\k$-algebra which
is concentrated in degree zero as a complex. Conversely, any generalized
Azumaya $\k$-algebra whose homology is concentrated in degree zero is
multiplicatively quasi-isomorphic to an ordinary Azumaya algebra.
\end{prop}
\begin{proof}
The first statement of the proposition is immediate. For the
second statement let $A$ be a generalized Azumaya $\k$-algebra
such that $H_*(A)=H_0(A)$. Then clearly $A$ is equivalent to
$H_0(A)$ in the homotopy category of differential graded algebras
over $\k$. Since $A$ is perfect as a $\k$-module it follows that
$H_0(A)$ is a finitely generated projective $\k$-module. The only
thing left to prove is that $H_0(A)$ is a faithful $\k$-algebra,
that is the unit map $\k\lra H_0(A)$ has zero kernel. Denote this
kernel by $I$, so there is a short exact sequence
\[
0\ra I\lra\k\lra H_0(A)\ra 0.
\]
Tensoring this short exact sequence with $H_0(A)$ over $\k$ and
taking into account that $H_0(A)$ is $\k$-projective we obtain the
following short exact sequence:
\[
0\ra I\otimes_\k H_0(A)\lra H_0(A)\lra H_0(A)\otimes_\k H_0(A)\ra
0.
\]
Since $H_0(A)$ is a ring the map $ H_0(A)\lra H_0(A)\otimes_\k H_0(A)$
is split by the multiplication map and it follows that
$I\otimes_\k H_0(A)=0$. In other words $I$ is $H_0(A)$-acyclic, or
equivalently $A$-acyclic. However, $\k$ is $A$-local by assumption,
which means that there is no nontrivial map into $\k$ from an $A$-acyclic
$\k$-module. Therefore $I=0$ and we are done.
\end{proof}

In the topological context it seems natural to relax the notion of
the Azumaya algebra somewhat.
\begin{defi}\label{defi:WeakTopAzumaya}
A \emph{weak topological Azumaya $R$-algebra}, or just a \emph{weak
Azumaya $R$-algebra}, is an $R$-algebra $A$ satisfying conditions (1)
and (2) of Definition~\ref{defi:AzumayaAlg} and such that the canonical
map $R\lra\THH_R(A,A)$ is the $A$-localization map of $R$ as an $R$-module.
\end{defi}

In the algebraic situation the analogue of the weak Azumaya $\k$-algebra
would be an algebra $A$ which is a finitely generated projective $\k$-module
and such that $A$ is an Azumaya algebra over $\k/I$ where, as above, $I$
is the kernel of the unit map $\k\lra A$.

Now we have the analogue of
Proposition~\ref{prop:AzumayaAlg-Char}.
\begin{prop}\label{prop:AzumayaAlg-weak}
An $R$-algebra $A$ is a weak Azumaya algebra if and only if the
following two conditions are satisfied.
\begin{enumerate}
\item
$A$ is a dualizable $R$-module.
\item
$\mu$ is a weak equivalence.
\end{enumerate}
\end{prop}
\proof
The proof is almost identical to that of
Proposition~\ref{prop:AzumayaAlg-Char}, the only new ingredient
being the weak equivalence
$$
\F_{R_A}(A,A)\cong\F_{R_A\wedge A}(A\wedge A,A)\simeq\F_A(A\wedge A,A)
\cong\F(A,A).
\eqno{\qed}
$$

\begin{rem}\label{rem:WeakTopAzumaya}
Notice that a weak topological Azumaya $R$-algebra $A$ is a
topological Azumaya algebra over $R_A$.
\end{rem}

A large supply of weak Azumaya algebras is provided by taking
endomorphisms of dualizable $R$-modules.
\begin{prop}\label{prop:Dualisable-Azumaya}
Let $E$ be a dualizable $q$-cofibrant $R$-module and $A=\F(E,E)$.
Then $A$ is a weak Azumaya algebra.
\end{prop}
\proof
It is clear that $A$ is a dualizable $R$-module. Furthermore $A$
acts canonically on the left on $E$ and therefore $A$ acts
\emph{on the right} on $E^*=F(E,R)$. Denote by $\tilde A$ the
cofibrant replacement of $A$ as an $R$-algebra and by
$\tilde{E^*}$ the cofibrant replacement of the right $\tilde
A$-module $E^*$. There is a map of $R$-algebras
\begin{equation}\label{eqn:op}
\tilde{A}^\op\lra\F(\tilde{E^*},\tilde{E^*})
\end{equation}
which is a weak equivalence since $E$ and hence $\tilde{E^*}$
are dualizable $R$-modules. Now the result follows from the
following chain of equivalences in the homotopy category of
$R$-algebras
$$
\tilde{A}\wedge\tilde{A}^\op\simeq\F(E,E)\wedge\F(\tilde{E^*},\tilde{E^*})
\simeq \F(E\wedge\tilde{E^*},E\wedge\tilde{E^*})\simeq\F(\tilde{A},\tilde{A}).
\eqno{\qed}
$$

We now give a rather general method of constructing nontrivial examples of
topological Azumaya algebras. Let $R$ be a commutative $\mathbb{S}$-algebra
whose coefficient ring $R_*$ is concentrated in even degrees and suppose that
$x\in R_d$ is not a zero divisor, so $d$ is even. Following Strickland~\cite{Str},
we consider the homotopy fibre sequence of $R$-modules
\begin{equation}\label{eqn:fib}
\Sigma^dR\xrightarrow{x}R\xrightarrow{\rho}R/x\xrightarrow{\beta}\Sigma^{d+1}R.
\end{equation}
Smashing this on the left with $R/x$ we obtain a homotopy fibre
sequence
\[
R/x\wedge\Sigma^dR\xrightarrow{x}R/x\wedge
R\xrightarrow{1\wedge\rho}R/x\wedge R/x.
\]
Since the first map is null this fibre sequence is split. By~\cite[lemma 3.7]{Str},
any such splitting determines a unital and associative product
$\phi\:R/x\wedge R/x\lra R/x$, and any other product $\phi'\:R/x\wedge R/x\lra R/x$
is of the form
\[
\phi'=\phi+u\circ(\beta\wedge\beta)
\]
for a unique element $u\in\pi_{2d+2}R/x$. If $\tau\:(R/x)^{\wedge 2}\lra(R/x)^{\wedge2}$
is the switch map, then $\phi\circ\tau$ is the multiplication opposite to $\phi$,
so
\begin{equation}\label{eqn:OppMult}
\phi\circ\tau=\phi+v\circ(\beta\wedge\beta)
\end{equation}
for a unique element $v\in\pi_{2d+2}R/x$.

Using the K\"unneth and universal coefficient spectral sequences,
it is straightforward to verify that as $R_*/x$-algebras,
\[
\pi_*R/x\wedge R/x^\op=\Lambda_{R_*/x}(\alpha), \quad
\pi_*\F(R/x,R/x)=\Lambda_{R/x_*}(\bar{\beta}),
\]
where $|\alpha|=d+1$ and $|\bar{\beta}|=-d-1$. We need to make an
explicit choice of generators $\alpha$ and $\bar{\beta}$. For
$\alpha$, consider the cofibre sequence
\[
R\wedge R/x^\op\xrightarrow{\rho\wedge1}R/x\wedge R/x^\op
\xrightarrow{\rho\wedge1}\Sigma^{d+1}R\wedge R/x^\op.
\]
On applying $\pi_*$ we obtain a short exact sequence
\[
0=\pi_{d+1}R\wedge R/x^\op\lra\pi_{d+1}R/x\wedge R/x^\op
\lra\pi_{d+1}\Sigma^{d+1}R\wedge R/x^\op=\pi_0R/x\ra0
\]
in which $\rho\in\pi_0R/x=\pi_{d+1}\Sigma^{d+1}R\wedge R/x^\op$
pulls back to a unique element $\alpha\in\pi_{d+1}R/x\wedge
R/x^\op$. We also take
\[
\bar{\beta}=\rho\circ\beta\:R/x\lra\Sigma^{d+1}R/x.
\]

Since $R/x$ is an $R$-ring spectrum, Lemma~\ref{lem:ActionMap}
shows that $\mu\:R/x\wedge R/x^\op\lra\F(R/x,R/x)$ is a map of
$R$-ring spectra.
\begin{prop}\label{prop:im}
We have
\[
\mu_*\alpha=v\bar{\beta}\in\pi_{d+1}\F(R/x,R/x).
\]
\end{prop}
\begin{proof}
Since $\phi\circ(\rho\wedge\rho)\:R\wedge R=R\lra R/x$ is
the unit map, we have
\[
\phi\circ\tau=\phi+v\circ(\beta\wedge\beta)
=\phi+v\phi\circ(\rho\wedge\rho)\circ(\beta\wedge\beta)
=\phi+v\phi\circ\bar{\beta}\wedge\bar{\beta}.
\]
Now consider the map
\[
\bar{\beta}'= \id\wedge\bar{\beta}\:R/x\wedge
R/x\lra\Sigma^{d+1}R/x\wedge R/x
\]
and recall that there are non-canonical equivalences of
$R$-modules
\[
R/x\wedge R/x\simeq R/x\vee\Sigma^{d+1}R/x, \quad
\Sigma^{d+1}R/x\wedge R/x\simeq\Sigma^{d+1}R/x\wedge R/x.
\]
Applying $\pi_*$ we obtain an induced map
\[
\bar{\beta}'_*\:\pi_*R/x\wedge R/x\lra\pi_*\Sigma^{d+1}R/x\wedge
R/x
\]
whose effect we wish to analyze in terms of the associated direct
sum decomposition
\[
\pi_*R/x\oplus\pi_*\Sigma^{d+1}R/x\lra
\pi_*\Sigma^{d+1}R/x\oplus\pi_*\Sigma^{2d+2}R/x.
\]
The natural $R_*$-module generating sets for the domain and
codomain are $1,\alpha$ and $\Sigma^{d+1}1,\Sigma^{2d+2}\alpha$.
The following result is independent of the actual choice of
splitting made above.
\begin{lem}\label{lem:bok}
With respect to the above decomposition, $\bar{\beta}'_*$
satisfies
\[
\bar{\beta}'_*(1)=0,\quad\bar{\beta}'_*(\alpha)=\Sigma^{d+1}1.
\]
\end{lem}
\proof
Smashing the fibre sequence of~\eqref{eqn:fib} with $R/x$ on the
left and taking into account that the map $x\:\Sigma^{d+1}R/x\lra
R/x$ is homotopic to zero we find that there is a split homotopy
fibre sequence
\[
R/x=R\wedge R/x\xrightarrow{\rho\wedge\id}R/x\wedge R/x
\xrightarrow{\beta\wedge\id}\Sigma^{d+1}R\wedge R/x=\Sigma^{d+1}R/x.
\]
Now the statement of the lemma becomes obvious after we take the composition
of $\beta\wedge\id$ with the split inclusion
$$
\Sigma^{d+1}R\wedge R/x\xrightarrow{\rho\wedge\id}\Sigma^{d+1}R/x\wedge R/x.
\eqno{\qed}
$$

Using~\eqref{eqn:OppMult} we find that
\begin{align*}
f=&\phi\circ(\id\wedge\phi)+
\phi\circ(v\wedge\phi)\circ(\id\wedge\bar{\beta}\wedge\bar{\beta})
\\
=&\phi\circ(\id\wedge\phi)+\phi\circ(v\bar{\beta}_*\wedge\bar{\beta}).
\end{align*}
Furthermore, the image of $\alpha\in\pi_{d+1}(R/x\wedge R/x^\op)$ in
$\pi_{-d-1}\F(R/x,R/x)=[R/x,\Sigma^{d+1}R/x]$ agrees with the composition
\[
\gamma\:R/x=R\wedge R/x\xrightarrow{\alpha\wedge\id}
\Sigma^{d+1}R/x\wedge R/x^\op\wedge
R/x\xrightarrow{f}\Sigma^{d+1}R/x.
\]
The composite
\[
\Sigma^{d+1}R\xrightarrow{\alpha}R/x\wedge
R/x^\op\xrightarrow{\phi}R/x
\]
is trivial for dimensional reasons, so to compute $\gamma$
only the component $\phi\circ(v\bar{\beta}'\wedge\bar{\beta})$
of $f$ is required. It follows from Lemma~\ref{lem:bok} that
the composition
\[
S^{d+1}\lra\Sigma^{d+1}R\xrightarrow{\alpha}
\Sigma^{d+1}R/x\wedge R/x^\op\xrightarrow{\bar{\beta}'}
R/x\wedge R/x^\op\xrightarrow{\phi}\Sigma^{d+1}R/x
\]
is the $(d+1)$-st suspension of the unit map, therefore
$\gamma=v\bar{\beta}$.
\end{proof}
\begin{cor}\label{cor:sd}
Suppose that the $R$-ring spectrum $R/x$ is actually an $R$-algebra
and the element $v\in\pi_{2d+2}R/x$ is invertible. Then $R/x$ is a
weak Azumaya $R$-algebra.
\end{cor}
\begin{proof}
By Lemma~\ref{lem:ActionMap}, $\mu\:R/x\wedge R/x^\op\lra\F(R/x,R/x)$
is a map of $R$-algebras, while Proposition~\ref{prop:im} implies that
the canonical $R$-algebra map $\mu$ is a weak equivalence. Since $R/x$
is a finite cell $R$-module, the result follows.
\end{proof}
\begin{cor}\label{cor:mod2}
$KU/2$ is a weak Azumaya $KU$-algebra and $v_1^{-1}MU/2$ is a weak
Azumaya $MU$-algebra.
\end{cor}
\begin{proof}
It was proved in~\cite{Laz} that the spectrum $KU/2$ is a $KU$-algebra
and $MU/2$ is an $MU$-algebra. Therefore our claim follows from
Corollary~\ref{cor:sd} and the formula for the commutators in $MU/2$
and $KU/2$ given in~\cite[theorem~2.12]{Str}.
\end{proof}
\begin{rem}\label{rem:nontrivial}
It seems likely that a weak Azumaya $R$ algebra of the form $R/x$ is
never isomorphic to an endomorphism algebra of a dualizable $R$-module.
Indeed, let us consider the case when the ring $R_*/x$ is a graded field
(this holds for $R=KU$ and $R/x=KU/2$). If $R/x$ were isomorphic to
$\F_R(M,M)$, where the $R$-module $M$ is dualizable, then being an
$R/x$-module, $M$ would split as a wedge of suspensions of $R/x$. But
then the hypercohomology spectral sequence
\[
\mathrm{E}_2^{*\,*}=\Ext_{R_*}^{*\,*}(M_*,M_*)\Lra\pi_*\F_R(M,M)
\]
clearly collapses and it follows that the graded $R_*$-module $\pi_*\F_R(M,M)$
contains odd degree elements. However $\pi_*\F_R(M,M)=R_*/x$, giving a
contradiction.
\end{rem}

The following conjecture is supported by the algebraic analysis of $A_\infty$
structures on quotients $R/x$ due to the second author~\cite{Laz2}.
\begin{con}\label{con:ext}
Let $R$ be a commutative $\mathbb{S}$-algebra such that $R_*$ is an evenly
graded ring and $x\in R_*$ is a nonzero divisor. Then any homotopy associative
$R$-ring spectrum structure on $R/x$ can be extended to an $R$-algebra structure.
\end{con}

Suppose provisionally that the above conjecture is true. Let $R$ be a commutative
$\mathbb{S}$-algebra whose homotopy is concentrated in even degrees with elements
$x\in R_d$, $v\in R_{2d+2}$ for which $x,v$ is a regular sequence in $R_*$. In
addition we suppose that~$2$ is invertible in $R_*$. Then by~\cite{Str}, $R/x$
has a unique commutative product $\phi\:R/x\wedge R/x\lra R/x$. Consider a new
product on $R/x$, say
\[
\phi'=\phi+v/2\circ(\beta\wedge\beta).
\]
From~\cite[Lemma~3.11]{Str} it follows that
\[
\phi'\circ\tau=\phi+v\circ(\beta\wedge\beta).
\]
By Conjecture~\ref{con:ext}, the homotopy multiplication $\phi'$ can be extended
to an $A_\infty$ structure. This leads to
\begin{cor}\label{cor:v^(-1)R/x-Azumaya}
If \emph{Conjecture~\ref{con:ext}} is true, then the $R$-algebra $v^{-1}R/x$ has
a structure of a topological Azumaya $R$-algebra.
\end{cor}

\section{Topological Hochschild cohomology of topological $K$-theory modulo~$2$}
\label{sec:THH-K(1)}

This section contains our main application of the Morita equivalence which
is the computation of the Hochschild cohomology of $K(1)$ at the prime~$2$.
We remark that McClure and Staffeldt~\cite{McCS} contain related results
obtained by different methods.

Let $\widehat{KU}_2$ denote the $2$-adic completion of the complex $K$-theory
spectrum with homotopy
\[
\pi_*\widehat{KU}_2=\Z_2[v_1,v_1^{-1}].
\]
Since we work in the categories $\mathscr{M}_{\widehat{KU}_2}$ and
$\mathscr{D}_{\widehat{KU}_2}$ and also $\mathscr{M}_{\mathbb{S}}$ and
$\mathscr{D}_{\mathbb{S}}$, we include subscripts indicating when smash
products, function spectra and topological Hochschild cohomology are
taken over $\widehat{KU}_2$ and omit them when they are taken over
$\mathbb{S}$; thus for example, $\THH=\THH_{\mathbb{S}}$.
\begin{theorem}\label{prop:THH-K(1)}
At the prime $2$, there is a weak equivalence of spectra
\[
\THH(K(1),K(1))\simeq\widehat{KU}_2.
\]
\end{theorem}
\begin{proof}
Observe that Corollary~\ref{cor:mod2} together with the fact that
$K(1)=\widehat{KU}_2/2$ implies that $K(1)$ is a topological
Azumaya $\widehat{KU}_2$-algebra and therefore
\[
\THH_{\widehat{KU}_2}(K(1),K(1))\simeq\widehat{KU}_2.
\]
Next consider the map
\[
l\:K(1)\wedge K(1)^\op\lra K(1)\wedge_{\widehat{KU}_2}K(1)^\op
\]
induced by the unit map $S\lra\widehat{KU}_2$. From~\cite{Nas} we
have
\[
K(1)_*K(1)^\op=
K(1)_*[a_0,t_k:k\geq1]/(a_0^2-t_1,v_1t_k^{2}-v_1^{2^k}t_k:k\geq1)
\]
as $K(1)_*$-algebras. Furthermore,
\[
K(1)^{\widehat{KU}_2}_*K(1)^\op=\pi_*K(1)\wedge_{\widehat{KU}_2}K(1)^\op
=\Lambda_{K(1)_*}(\tau_0).
\]
\begin{lem}\label{lem:THH-K(1)}
The map $l_*\:K(1)_*K(1)^\op\lra K(1)_*^{\widehat{KU}_2}K(1)^\op$
is surjective with kernel generated by the elements $t_k$.
\end{lem}
\begin{proof}
First we prove the surjectivity of $l_*$. If $l_*$ is not
surjective then its image is the subring $K(1)_*$ of
$K(1)_*^{\widehat{KU}_2}K(1)^\op$ and the dual map
\begin{equation}\label{eqn:dual}
\delta\:\Hom_{K(1)_*}(K(1)_*^{\widehat{KU}_2}K(1)^\op,K(1)_*)\lra
\Hom_{K(1)_*}(K(1)_* K(1)^\op,K(1)_*)
\end{equation}
must factor through $\Hom_{K(1)_*}(K(1)_*,K(1)_*)=K(1)_*$. Notice
that the map of~\eqref{eqn:dual} agrees with the forgetful map
\[
[K(1)^\op,K(1)]_*^{\widehat{KU}_2}\cong[K(1),K(1)]_*^{\widehat{KU}_2}
\lra[K(1)^\op,K(1)]_*=[K(1),K(1)]_*.
\]
But clearly the mod~$2$ Bockstein operation is a self-map of $K(1)$
considered either as an $\mathbb{S}$-module or an $\widehat{KU}_2$-module.
This shows that the map $\delta$ of~\eqref{eqn:dual} cannot factor through
$K(1)$ and so $l_*$ is onto.

Next we check that $l_*(t_k)=0$. This is equivalent to asserting that
the map $\delta$ of~\eqref{eqn:dual} factors through
\[
\Hom_{K(1)_*}(K(1)_*K(1)^{\op}\otimes_{\Sigma_*}K(1)_*,K(1)_*),
\]
where $\Sigma_*$ is the $K(1)_*$-subalgebra in $K(1)_*K(1)^\op$ generated
by the $t_k$'s. Since the image of $l_*$ is generated by the mod~$2$ Bockstein
operation $Q_0$, it suffices to show that $Q_0(t_k)=0$ for all $k$. Since
$Q_0$ is a derivation it follows from the well-known fact that the algebra
$\Sigma_*$ is semisimple.

Therefore $\ker l_*$ is an ideal containing all of $t_k$'s. It cannot contain
$a_0$ since in that case $l_*$ would not be onto. The only remaining possibility
is that $\ker l_*$ is generated by $t_k$'s.
\end{proof}

Now consider two spectral sequences
\begin{align}
\mathrm{E}(1)^{*\,*}_2=\Ext^{*\,*}_{K(1)_*K(1)}(K(1)_*,K(1)_*)
&\Lra\THH^{*}(K(1),K(1)), \label{eqn:SS-first}\\
\mathrm{E}(2)^{*\,*}_2=\Ext^{*\,*}_{K(1)^{\widehat{KU}_2}_*K(1)}(K(1)_*,K(1)_*)
&\Lra\THH^*_{\widehat{KU}_2}(K(1),K(1)). \label{eqn:SS-second}
\end{align}
It is now immediate that the $\mathrm{E}_2$-term of~\eqref{eqn:SS-second}
is
\[
\mathrm{E}(2)^{*\,*}_2=
\Ext^{*\,*}_{\Lambda_{K(1)_*}(\tau_0)}(K(1)_*,K(1)_*)=K(1)_*[[y]],
\]
where $y$ has bidegree $(1,1)$. Once again using the fact that $\Sigma_*$
is semisimple, we obtain the following base change isomorphism
\begin{align*}
\mathrm{E}(1)^{*\,*}_2
&\iso\Ext^{*\,*}_{K(1)_*K(1)\otimes_{\Sigma_*}K(1)_*}
(K(1)_*\otimes_{\Sigma_*}K(1)_*,K(1)_*) \\
&=\Ext^{*\,*}_{\Lambda_{K(1)_*}(a_0)}(K(1)_*,K(1)_*) \\
&=K(1)_*[[y]].
\end{align*}
It follows that the map $l_*$ induces an isomorphism between the
spectral sequences of~\eqref{eqn:SS-first}
and~\eqref{eqn:SS-second}. Therefore
\[
\THH(K(1),K(1))\simeq\THH_{\widehat{KU}_2}(K(1),K(1))\simeq\hat{KU}_2,
\]
giving the result.
\end{proof}

It is natural to ask what happens when $p$ is odd. First of all, the answer
depends on whether we consider the $2$-periodic or the $2(p-1)$-periodic
version of $K(1)$. From~\cite{AB&BR:E-infinity} it is known that the Adams
summand $E(1)$ of $KU_{(p)}$ and its $p$-completion $\hat{E(1)}$ have unique
commutative $\mathbb{S}$-algebra structures. The $2(p-1)$-periodic $K(1)$-theory
can be shown to be a homotopy commutative $E(1)$-ring spectrum using the
methods of~\cite{Wurgler:K(n),EKMM,Str}. Since the action map
$\mu\:K(1)\wedge_{E(1)}K(1)\lra\mathrm{F}_{E(1)}(K(1),K(1))$ has a homotopy
factorization through the augmentation map to $K(1)$, it is not an equivalence
and therefore $K(1)$ cannot be a topological Azumaya algebra. On the other hand,
if we considers the $2$-periodic theory then there exists a homotopy associative
multiplication on $K(1)$ for which the canonical map
$K(1)\wedge_{KU}K(1)^\op\lra\F_{KU}(K(1),K(1))$ is a weak equivalence.
If Conjecture~\ref{con:ext} were true then in this case we would also have
$\THH(K(1),K(1))\simeq\hat{KU}_2$. To summarize, it seems likely that for
an odd prime~$p$ the result depends heavily on what $\mathbb{S}$-algebra
structure is chosen on $K(1)$. This naturally leads to
\begin{prob}\label{prob:structures}
For an odd prime $p$, classify non-isomorphic $\mathbb{S}$-algebra or
$\widehat{KU}_p$-algebra structures on the spectrum $K(1)$.
\end{prob}

To illustrate the subtlety of this problem consider again the case
$p=2$. We claim that as an $\widehat{KU}_2$-algebra $K(1)$ is not
isomorphic to its opposite $K(1)^\op$. To see this, consider the
forgetful map
\begin{equation}\label{eqn:forget}
\mu\:K(1)_*K(1)\lra K(1)_*^{\widehat{KU}_2}K(1).
\end{equation}
By lemma~8 of~\cite{Nas}, we know that
$K(1)_*K(1)=\Sigma_*[a_0]/(a_0^2-t_1-v_1)$ and arguments similar to
the ones used in the proof of Theorem~\ref{prop:THH-K(1)} show that
$K(1)_*^{\widehat{KU}_2}K(1)\cong K(1)_*[\tau_0]/(\tau_0^2-v_1)$ and
the map $\mu$ of~\eqref{eqn:forget} is onto, with kernel generated
by the elements $t_i$.

If $K(1)$ were isomorphic to $K(1)^\op$ as an $\widehat{KU}_2$-algebra,
then the $K(1)_*$-algebras $K(1)_*^{\widehat{KU}_2}K(1)$ and
$K(1)_*^{\widehat{KU}_2}K(1)^\op$ would be isomorphic. But this is
clearly false since $K(1)_*^{\widehat{KU}_2}K(1)^\op$ contains an
exterior generator while $K(1)_*^{\widehat{KU}_2}K(1)$ does not.

In effect the above argument shows that $K(1)$ and $K(1)^\op$ are not
isomorphic even as $\widehat{KU}_2$-ring spectra. On the other hand,
by~\cite[lemma~9]{Nas}, the ring spectrum $K(1)$ possesses a nontrivial
anti-involution from which it follows that $K(1)$ is isomorphic to
$K(1)^\op$ as a ring spectrum. We are led to make the following conjecture.
\begin{con}\label{con:K(1)-Salg}
For $p=2$, there are precisely two structures of an $\widehat{KU}_2$-algebra
on $K(1)$ and these are opposite to each other. Moreover, the spectrum $K(1)$
possesses a unique structure of an $\mathbb{S}$-algebra.
\end{con}

Note that this
conjecture implies that the $\widehat{KU}_2$-algebra
structure on $K(1)\wedge_{\widehat{KU}_2}K(1)^\op$ does not depend on
the choice of the $\widehat{KU}_2$-algebra structure on $K(1)$. Thus
$\THH_{\widehat{KU}_2}(K(1),K(1))$ is also independent of this choice.
This is in agreement with Theorem~\ref{prop:THH-K(1)}.

\section{Future directions}\label{sec:Future}

Perhaps one of the most interesting ramifications of our constructions
lies in the theory of the \emph{Brauer group}. The classical Brauer
group of a ring $R$ is formed by classes of Morita equivalent Azumaya
algebras over~$R$. Thus, endomorphism rings of finitely generated
projective modules are considered trivial in the Brauer group. Our
version of Morita theory suggests considering the topological Brauer
group of an $\mathbb{S}$-algebra $R$ as generated by topological Azumaya
$R$-algebras modulo the endomorphism algebras of dualizable $R$-modules
(this yields a set since equivalence classes of retracts of finite cell
$R$-modules clearly form a set and so equivalence classes of algebra
structures on them also do). Note that for any finite cell $R$-modules
$A,B$, there is a weak equivalence
\[
\mathrm{F}_R(A,A)\wedge\mathrm{F}_R(B,B)\simeq\mathrm{F}_R(A\wedge B,A\wedge B).
\]
This implies that the set of endomorphism algebras of finite cell
$R$-modules is closed with respect to smash products and so it
forms a submonoid inside the monoid of all Azumaya algebras.

In this sense, the $\widehat{KU}_2$-algebra $KU/2$ gives a nontrivial
example of a topological Azumaya algebra, see Remark~\ref{rem:nontrivial}.
We intend to return to this topic in future work; it seems likely there
are connections with the Galois theory of commutative $\mathbb{S}$-algebras
initiated by John Rognes, see~\cite{AB&BR:Galois}.

We note that all of these constructions also make sense in the category
of differential graded modules over a fixed commutative ring. This extends
the notion of the Brauer group from the purely algebraic context. There
are examples of nontrivial generalized Azumaya algebras which are \emph{not}
Azumaya algebras in the classical sense. We restrict ourselves with giving
one example; a more detailed discussion can be found in~\cite{Laz2}.
This example provides an algebraic analogue of the $\widehat{KU}_p$-algebra
$KU/p$; however, in distinction to the topological case, it does not occur
for the prime~$2$.
\begin{example}\label{ex:Azumaya}
Let $R=\Z[v,v^{-1]}$ be the ring of Laurent polynomials in one variable~$v$
of degree~$2$ and let $p$ be an odd prime. Consider the differential graded
algebra $A$ over $R$ with a single generator $y$ in degree~$1$ subject to
the relation $y^2=v$ and whose differential given by $dy=p$. Then
\[
A\,{\mathop{\otimes}^{\mathrm{L}}}_R\,A^\op\simeq\RHom_R(A,A),
\]
and furthermore $A$ is clearly dualizable as an $R$-module and therefore
$A$ is a (weak) Azumaya $R$-algebra.
\end{example}

\Addresses\recd
\end{document}